\documentclass[11pt,reqno]{amsart}

\usepackage{epsfig}
\usepackage{graphicx}
\usepackage{mathrsfs}
\usepackage{subfigure}
\usepackage{verbatim}
\usepackage{amssymb}
\usepackage{latexsym}
\usepackage{amsmath}
\usepackage{amsfonts}
\usepackage{amsbsy}
\usepackage{amscd}
\usepackage{stmaryrd}
\usepackage[dvipsnames,usenames]{color} 
\usepackage{color}
\usepackage{epstopdf}
\usepackage{hyperref}
\usepackage{cite}

\newcommand{\R}{{\mathbb R}}

\newtheorem{prop}{Proposition}[section]
\newtheorem{lem}[prop]{Lemma}

\newtheorem{defi}[prop]{Definition}

\newtheorem{thm}[prop]{Theorem}

\usepackage[misc,geometry]{ifsym}

\begin{document}
\baselineskip=16pt

\title[A nonlinear wave equation on graphs]{Application of Rothe's method to a nonlinear wave equation on graphs}

\author[Y. Lin]{Yong Lin}
\author[Y. Xie]{Yuanyuan Xie~\Letter}



\begin{abstract}
We study a nonlinear wave equation on finite connected weighted graphs. Using Rothe's and energy methods, we prove the existence and uniqueness of solution under certain assumption. For linear wave equation on graphs, Lin and Xie \cite{Lin-Xie} obtained the existence and uniqueness of solution. The main novelty of this paper is that the wave equation we considered has the nonlinear damping term $|u_t|^{p-1}\cdot u_t$ ($p>1$).
\end{abstract}
\maketitle

\section{Introduction}\label{S:intr}
\setcounter{equation}{0}

A graph is an ordered pair $(V,E)$ with $V$ being a set of vertices and $E$ being a set of edges. Let $\mu: V\to (0,\infty)$ be the vertex measure. Also, let $\omega: V\times V\to (0,\infty)$ be the edge weight function satisfying positivity and symmetry, that is, $\omega_{xy}>0$ and $\omega_{xy}=\omega_{yx}$ for any $xy\in E$. We write $y\sim x$ if $xy\in E$.
Define
\begin{equation*}
D_\mu:=\max\Big\{\frac{1}{\mu(x)}\sum_{y\sim x}\omega_{xy}: x\in V\Big\}.
\end{equation*}
The quadruple $G=(V,E,\mu,\omega)$ will be referred as a weighted graph.
In this paper, the graphs we consider are finite connected weighted.

Let $C(V):=\{v: V\to \R\}$. Define the $\mu$-Laplacian $\Delta$ of $v\in C(V)$ by
\begin{equation*}\label{e:mu_Laplacian}
\Delta v(x)=\frac{1}{\mu(x)}\sum_{y\sim x}\omega_{xy}\big(v(y)-v(x)\big).
\end{equation*}
We denote the associated gradient form by
\begin{equation*}\label{e:Gamma}
\Gamma(v_1,v_2)(x)=\frac{1}{2\mu(x)}\sum_{y\sim x}\omega_{xy}\big(v_1(y)-v_1(x)\big)\big(v_2(y)-v_2(x)\big).
\end{equation*}
Let $|\nabla v|^2(x):=\Gamma(v,v)(x)$, and $|\nabla v|(x)$ be the length of $\Gamma$. Also, write
\begin{equation*}\label{e:integration_of_function}
\int_V v\,d\mu=\sum_{x\in V}\mu(x)v(x)\qquad\mbox{ for any }v\in C(V).
\end{equation*}

For any non-empty domain $\Omega\subseteq V$, let
\begin{equation*}\label{e:defi_boundary_and_interior_of_Omega}
\partial\Omega:=\{y\in\Omega: \mbox{there exists }x\in V\backslash\Omega\mbox{ such that }xy\in E\}
\quad\mbox{and}\quad\Omega^\circ:=\Omega\setminus\partial\Omega.
\end{equation*}
For any real function $v$ on $\Omega^\circ$, we extend $v$ to $V$ by letting $v(x)=0$ for any $x\in V\backslash\Omega^\circ$. Set $\Delta_\Omega v=(\Delta v)|_{\Omega^\circ}$, we call $\Delta_\Omega$
the \textit{Dirichlet Laplacian} on $\Omega^\circ$. Then
\begin{equation*}\label{defi:Dirichlet_Lap}
\Delta_\Omega v(x)=\frac{1}{\mu(x)}\sum_{y\sim x} \omega_{xy}\big(v(y)-v(x)\big) \qquad\mbox{ on }\Omega^\circ,
\end{equation*}
where $v$ vanishes on $V\backslash\Omega^\circ$. Clearly, the operator $-\Delta_\Omega$ is positive and self-adjoint (see \cite{Grigoryan_2018, Weber_2012}).

Let $p>1$ be a constant. For give functions $f:[0,\infty)\times\Omega^\circ\to \R$, and $g, h: \Omega^\circ\to\R$, we study the problem
\begin{equation}\label{e:wave_eq_on_graph}
\left\{
\begin{aligned}
&u_{tt}-\Delta_\Omega u+|u_t|^{p-1}\cdot u_t=f,\quad&&t\ge 0, x\in\Omega^\circ,\\
&u|_{t=0}=g,\quad&& x\in\Omega^\circ,\\
&u_t|_{t=0}=h,\quad&&x\in\Omega^\circ,\\
&u=0,\quad&& t\ge 0, x \in\partial\Omega,
\end{aligned}
\right.
\end{equation}
where $f$ is continuous with respect to $t$.

\begin{defi}
We call $u=u(t,x)$ a {\rm solution} of \eqref{e:wave_eq_on_graph} on $[0,T]\times\Omega$ if $u$ is twice continuously differentiable with respect to $t$, and \eqref{e:wave_eq_on_graph} holds.
\end{defi}

The problem \eqref{e:wave_eq_on_graph} has been studied by Lions \cite{Lions_1969} who gave the existence and uniqueness of solution on $\R^d$. On metric graphs, Friedman and Tillich \cite{Friedman-Tillich_2004} studied the wave equation whose Laplacian is based on the edge. Recently, the authors \cite{Lin-Xie} considered the linear wave equation on graphs, and obtained the existence result of solution. The main difference between this paper and \cite{Lin-Xie} is that the problem \eqref{e:wave_eq_on_graph} has the nonlinear damping term $|u_t|^{p-1}\cdot u_t$. In this case, it is much harder to study the existence of solution.

In recent years, various partial differential equations have also been extensively studied on graphs. Using variational method, Grigoryan et al. \cite{Grigoryan-Lin-Yang_2016_1, Grigoryan-Lin-Yang_2016_2, Grigoryan-Lin-Yang_2017} gave existence results of the solution of Yamabe type equation, Kazdan-Warner equation and some nonlinear equations. Lin and Wu \cite{Lin-Wu_2017} considered a semilinear heat equation, and obtained the existence and nonexistence results of global solution. For more relevant results, please refer to \cite{Han-Shao-Zhao_2020, Huang-Lin-Yau} and their references.

In this paper, using Rothe's method that was originally introduced by Rothe \cite{Rothe_1930} for the study of parabolic equation, we obtain the solution of \eqref{e:wave_eq_on_graph} exists globally. After 1930, using this method, many authors
(e.g., \cite{Rektorys_1971, Kacur_1984}) obtained existence results for solutions to parabolic and hyperbolic equations.

Now, we briefly introduced Rothe's method. For any $T>0$, divide $[0,T]$ into $n$ equidistant subintervals $[t_{i-1}, t_i]$ with $t_0=0, t_n=T$ and $t_i=i {\delta}$ for $i\in \Lambda:=\{1,\ldots,n\}$. For $i\in\Lambda$, let $u_{n,0}, u_{n,-1}, f_{n,i}$ be defined as in subsection \ref{SS:some_priori_est}, and solve successively $n$ equations
\begin{equation*}
(u_{n,i}-2u_{n,i-1}+u_{n,i-2})/{{\delta}^2}-\Delta_\Omega u_{n,i}
+(u_{n,i}-u_{n,i-1})/{{\delta}}\cdot\big|(u_{n,i}-u_{n,i-1})/{{\delta}}\big|^{p-1}
=f_{n,i}\qquad\mbox{on }\Omega^\circ.
\end{equation*}
Using $\{u_{n,i}\}_{i\in\Lambda}$, we can construct Rothe's functions as following
\begin{equation*}
u^{(n)}(t,x)=u_{n,i-1}(x)+(t-t_i)\cdot(u_{n,i}(x)-u_{n,i-1}(x))/{\delta}
\qquad i\in\Lambda\mbox{ and }t\in[t_{i-1}, t_i].
\end{equation*}
Under certain assumption, we prove $\{u^{(n)}(t,x)\}$ converges to $u$, where $u$ is a solution of \eqref{e:wave_eq_on_graph}.

Throughout this paper, let $C_{\Omega^\circ}:=C(\Omega^\circ)>0$ be a constant depending only on $\Omega^\circ$.
Similarly, let $C_\Omega:=C(\Omega)>0$ and $C_{\Omega, p}:=C(\Omega, p)>0$.

Assume that for positive constants $\gamma$ and $C_{\Omega^\circ}$, the following holds
\begin{equation}\label{eq:condition_on_f}
\|f(s_1,\cdot)-f(s_2,\cdot)\|_{L^2(\Omega^\circ)}\le C_{\Omega^\circ}\cdot|s_1-s_2|^\gamma \quad\mbox{ for any }s_1, s_2\in[0,\infty).
\end{equation}
Now we state our main result.

\begin{thm}\label{T:sol_of_wave_equ}
Let $G=(V, E,\mu,\omega)$ be a finite connected weighted graph, and let $\Omega\subseteq V$ be a domain satisfying $\Omega^\circ\neq \emptyset$. If \eqref{eq:condition_on_f} holds, then \eqref{e:wave_eq_on_graph} has a unique global solution.
\end{thm}

We introduce Green's formula and Sobolev embedding theorem in Section \ref{S:pre}.
Theorem \ref{T:sol_of_wave_equ} will be proved in Section \ref{S:proof_of_sol_of_wave_equ}.

\section{Preliminaries}\label{S:pre}
\setcounter{equation}{0}

Let $G=(V, E, \mu, \omega)$ be a finite connected weighted graph, and $\Omega\subseteq V$ be a domain such that $\Omega^\circ$ is non-empty.

\begin{lem}\label{L:Green_formula}(Green's formula)\cite{Grigoryan_2018}
For any real functions $w, v$ on $\Omega^\circ$, we have
\begin{eqnarray*}
\int_{\Omega^\circ}\Delta_\Omega w\cdot v\,d\mu=-\int_\Omega\Gamma(w,v)\,d\mu.
\end{eqnarray*}\label{L:Green_formula}
\end{lem}

For $q\in[1,\infty)$, let $L^q(\Omega)$ is a space of all real-valued functions on $V$ whose norm
$\|v\|_{L^q}:=\{\int_\Omega|v|^q\,d\mu\}^{1/q}$ is finite.
For $q=\infty$, denote
$$L^\infty(\Omega):=\big\{v\in C(V):\sup\limits_{x\in\Omega}|v(x)|<\infty\big\}.$$
with norm $\|v\|_{L^\infty(\Omega)}=\sup\limits_{x\in\Omega}|v(x)|$.
It is easy to see that $L^q(\Omega)$ is a Banach space. Moreover, $L^2(\Omega)$ is a Hilbert space with the following inner product
\begin{eqnarray*}\label{eq:inner_pro_L2}
(w,v)=\int_\Omega w(x) v(x)\,d\mu\quad\mbox{ for }w,v\in L^2(\Omega).
\end{eqnarray*}
Let
$$W^{1,2}(\Omega):=\{v\in C(V): \int_\Omega (|\nabla v|^2+|v|^2)\,d\mu<\infty\}$$
with norm
\begin{eqnarray}\label{e:defi_W_12}
\|v\|_{W^{1,2}(\Omega)}=\Big(\int_\Omega (|\nabla v|^2+|v|^2)\,d\mu \Big)^{1/2}.
\end{eqnarray}
Let $C_0(\Omega):=\{v\in C(\Omega): v=0 \mbox{ on }\partial\Omega\}$. We complete $C_0(\Omega)$ under the norm \eqref{e:defi_W_12} and denote the completed space by $W^{1,2}_0(\Omega)$. Clearly $W^{1,2}_0(\Omega)$ is a Hilbert space under inner product
\begin{equation*}\label{eq:inner_product_W}
(w,v)_{W_0^{1,2}(\Omega)}=\int_\Omega(\Gamma(w,v)+wv)\,d\mu\qquad\mbox{ for any }w,v\in W_0^{1,2}(\Omega).
\end{equation*}

Since $\Omega$ is finite, the dimension of $W_0^{1,2}(\Omega)$ is finite. A graph $G$ is said to be \textit{locally finite} if for any $x\in V$, $\#\{y\in V: xy\in E\}$ is finite. It is obvious that a finite graph is locally finite. So we state the Sobolev embedding theorem (see \cite[Theorem 7]{Grigoryan-Lin-Yang_2016_1}) for finite graph.

\begin{thm}\label{T:Sobolev_embedding_thm}
Let $(V, E)$ be a finite graph, and $\Omega\subseteq V$ be a domain satisfying $\Omega^\circ\neq \emptyset$.
Then $W_0^{1,2}(\Omega)\hookrightarrow L^q(\Omega)$ for all $q\in[1,\infty]$. Particularly, there exists constant $C_\Omega$ such that
\begin{eqnarray*}
\|v\|_{L^q(\Omega)}\le C_\Omega\|\nabla v\|_{L^2(\Omega)}\qquad\mbox{ for all }q\in[1,\infty]\mbox{ and all }v\in W_0^{1,2}(\Omega).
\end{eqnarray*}
Moreover, $W_0^{1,2}(\Omega)$ is precompact, that is, a bounded sequence in $W_0^{1,2}(\Omega)$ contains a convergent subsequence.
\end{thm}

\section{Proof of Theorem~\ref{T:sol_of_wave_equ}}\label{S:proof_of_sol_of_wave_equ}
\setcounter{equation}{0}

In this section, we show that there exists a unique global solution of \eqref{e:wave_eq_on_graph}.
In subsection \ref{SS:some_priori_est}, we set up some priori estimates that will be used in the proof of Theorem~\ref{T:sol_of_wave_equ}.

\subsection{Some priori estimates}\label{SS:some_priori_est}

For any $T>0$, let $\{t_i\}_{i=0}^n$ be an equidistant partition of times interval $[0,T]$ satisfying
$t_0=0$, $t_n=T$, and $t_i=i {\delta}$ for $i\in \Lambda:=\{1,\ldots, n\}$. Let
$$u_{n,0}(x):=g(x),\quad u_{n,-1}(x):=g(x)-{\delta} h(x),\quad
f_{n,i}(x):=f(t_i,x)\quad\mbox{ for }i\in\Lambda, x\in\Omega^\circ,
$$
and $u_{n,0}(x)=u_{n,-1}(x)=0$ on $\partial\Omega$.

For $p>1$, define the functional $\mathcal{J}_1$ from $W_0^{1,2}(\Omega)$ to $\R$ as
\begin{eqnarray*}
\begin{aligned}
\mathcal{J}_1(u)
=&\int_{\Omega^\circ}(u-4u_{n,0}+2u_{n,-1})/{{\delta}^2}\cdot u\,d\mu+\int_{\Omega}|\nabla u|^2\,d\mu\\
&+2{\delta}/(p+1)\cdot\int_{\Omega^\circ}|(u-u_{n,0})/{{\delta}}|^{p+1}\,d\mu-2\int_{\Omega^\circ} f_{n,1}\cdot u\,d\mu.
\end{aligned}
\end{eqnarray*}

\begin{lem}
$\mathcal{J}_1(u)$ attains its minimum $u_{n,1}\in W_0^{1,2}(\Omega)$, and $u_{n,1}$ is the unique solution of
\begin{eqnarray}\label{e:BVP_1}
(u-2u_{n,0}+u_{n,-1})/{{\delta}^2}-\Delta_\Omega u+|(u-u_{n,0})/{{\delta}}|^{p-1}\cdot (u-u_{n,0})/{{\delta}}=f_{n,1}
\quad\mbox{on }\Omega^\circ.
\end{eqnarray}
\end{lem}

\begin{proof}
This proof consists two parts.

{\bf Part 1} We show that $\mathcal{J}_1(u)$ attains its minimum $u_{n,1}\in W_0^{1,2}(\Omega)$.
Using H\"older inequality, we obtain
\begin{eqnarray*}
\begin{aligned}
\mathcal{J}_1(u)
\ge&\int_{\Omega}|\nabla u|^2\,d\mu+{2{\delta}}/{(p+1)}\cdot\int_{\Omega^\circ}|(u-u_{n,0})/{{\delta}}|^{p+1}\,d\mu
-\int_{\Omega^\circ}|(2u_{n,0}-u_{n,-1})/{\delta}+\delta f_{n,1}|^2\,d\mu\\
\ge&-\int_{\Omega^\circ}|g/{\delta}+ h+{\delta}\cdot f({\delta},x)|^2\,d\mu,
\end{aligned}
\end{eqnarray*}
and so $\mathcal{J}_1$ has a lower bound on $W_0^{1,2}(\Omega)$. Further, $\inf_{u\in W_0^{1,2}(\Omega)}\mathcal{J}_1$ is finite.

Taking a sequence of functions $\{u_k\}\subseteq W_0^{1,2}(\Omega)$ such that $\mathcal{J}_1(u_k)\to a_1:=\inf_{u\in W_0^{1,2}(\Omega)}\mathcal{J}_1$. That is, $|\mathcal{J}_1-a_1|<\epsilon_1$ for some $\epsilon_1>0$, and so
\begin{eqnarray*}
\int_{\Omega}|\nabla u_k|^2\,d\mu
\le\int_{\Omega^\circ}|g/{\delta}+h+{\delta} f({\delta},x)|^2\,d\mu+a_1+\epsilon_1,
\end{eqnarray*}
which, together with Theorem \ref{T:Sobolev_embedding_thm}, yields $u_k$ is bounded in $W_0^{1,2}(\Omega)$. Also,
there exist a function $u_{n,1}\in W_0^{1,2}(\Omega)$ and a subsequence $\{u_{k_j}\}$ such that $u_{k_j}\to u_{n,1}$ in $W_0^{1,2}(\Omega)$. Further, $\|u_{k_j}\|_{W^{1,2}(\Omega)}\to \|u_{n,1}\|_{W^{1,2}(\Omega)}$. Since
\begin{eqnarray*}
\big|\|u_{k_j}\|_{L^2(\Omega)}-\|u_{n,1}\|_{L^2(\Omega)}\big|\le \|u_{k_j}-u_{n,1}\|_{L^2(\Omega)}
\le \|u_{k_j}-u_{n,1}\|_{W^{1,2}(\Omega)},
\end{eqnarray*}
we obtain
\begin{eqnarray}\label{e:u_nk_to_u_infty_in_L2_W12}
\|u_{k_j}\|^2_{L^2(\Omega)}\to \|u_{n,1}\|^2_{L^2(\Omega)}\quad\mbox{and}\quad
\|\nabla u_{k_j}\|^2_{L^2(\Omega)}\to \|\nabla u_{n,1}\|^2_{L^2(\Omega)}.
\end{eqnarray}
Moreover, $u_{k_j}\to u_{n,1}$ on $\Omega$. Based on the above results, we get
\begin{eqnarray*}
\mathcal{J}_1(u_{n,1})=\lim_{j\to\infty}\mathcal{J}_1(u_{k_j})=a_1.
\end{eqnarray*}
This proves that $\mathcal{J}_1$ attains its minimum $u_{n,1}\in W_0^{1,2}(\Omega)$.

{\bf Part 2} We prove that $u_{n,1}$ is the unique solution of \eqref{e:BVP_1}. For any $\psi\in W_0^{1,2}(\Omega)$,
\begin{eqnarray*}
\begin{aligned}
0=&\frac{d}{d\eta}\Big|_{\eta=0}\mathcal{J}_1(u_{n,1}+\eta\psi)\\
=&2\int_{\Omega^\circ}\Big((u_{n,1}-2u_{n,0}+u_{n,-1})/{{\delta}^2}-\Delta_\Omega u_{n,1}\\
&\qquad+\big|(u_{n,1}-u_{n,0})/{{\delta}}\big|^{p-1}\cdot (u_{n,1}-u_{n,0})/{{\delta}}-f_{n,1}\Big)\cdot\psi\,d\mu.
\end{aligned}
\end{eqnarray*}
This proves $u_{n,1}$ is a solution of \eqref{e:BVP_1}.

Let $u_{n,1}$ and $\breve{u}$ be two solution of \eqref{e:BVP_1}. Then for $p>1$,
\begin{eqnarray}\label{e:equation_on_u1}
\begin{aligned}
&(u_{n,1}-\breve{u})/{{\delta}^2}-\Delta_\Omega(u_{n,1}-\breve{u})
+\big|(u_{n,1}-u_{n,0})/{{\delta}}\big|^{p-1}\cdot(u_{n,1}-u_{n,0})/{{\delta}}\\
-&\big|(\breve{u}-u_{n,0})/{{\delta}}\big|^{p-1}\cdot(\breve{u}-u_{n,0})/{{\delta}}=0\qquad\mbox{ on }\Omega^\circ.
\end{aligned}
\end{eqnarray}
For any $x_0\in \Omega^\circ$, if $(u_{n,1}-\breve{u})(x_0)\ge 0$, then $-\Delta(u_{n,1}-\breve{u})(x_0)\ge0$, and
\begin{eqnarray*}
\big|(u_{n,1}-u_{n,0})/{{\delta}}\big|^{p-1}\cdot(u_{n,1}-u_{n,0})/{{\delta}}
-\big|(\breve{u}-u_{n,0})/{{\delta}}\big|^{p-1}\cdot(\breve{u}-u_{n,0})/{{\delta}}\ge0.
\end{eqnarray*}
Combining these with \eqref{e:equation_on_u1}, we get $u_{n,1}(x_0)=\breve{u}(x_0)$.
Then $u_{n,1}=\breve{u}$ on $\Omega^\circ$ follows from that $x_0$ is arbitrary. This completes the proof.
\end{proof}

Successively, for $i\in \Lambda\backslash\{1\}$, consider the functionals $\mathcal{J}_i$ from $W_0^{1,2}(\Omega)$ to $\R$:
\begin{eqnarray*}
\begin{aligned}
\mathcal{J}_i(u)
=&\int_{\Omega^\circ}(u-4u_{n,i-1}+2u_{n,i-2})/{{\delta}^2}\cdot u\,d\mu+\int_{\Omega}|\nabla u|^2\,d\mu\\
&+2{\delta}/(p+1)\cdot\int_{\Omega^\circ}\big|(u-u_{n,i-1})/{{\delta}}\big|^{p+1}\,d\mu-2\int_{\Omega^\circ} f_{n,i}\cdot u\,d\mu.
\end{aligned}
\end{eqnarray*}
Similarly, $\mathcal{J}_i$ attains its minimum $u_{n,i}\in W_0^{1,2}(\Omega)$, and $u_{n,i}$ solves uniquely
\begin{eqnarray}\label{e:BVP_i}
(u-2u_{n,i-1}+u_{n,i-2})/{{\delta}^2}-\Delta_\Omega u+(u-u_{n,i-1})/{{\delta}}\cdot\big|(u-u_{n,i-1})/{{\delta}}\big|^{p-1}
=f_{n,i}\qquad\mbox{on }\Omega^\circ.
\end{eqnarray}
Let $u_{n,i}(x)$ be the approximation of $u(t,x)$, which is the solution of \eqref{e:wave_eq_on_graph}, at $t=t_i$.
We denote
\begin{eqnarray}\label{e:defi_delta_u_ni}
w_{n,i}(x):=(u_{n,i}(x)-u_{n,i-1}(x))/{\delta}\quad\mbox{ for }i\in\Lambda\cup\{0\},
\end{eqnarray}
\begin{eqnarray}\label{e:defi_delta2_u_ni}
z_{n,i}(x):=(w_{n,i}(x)-w_{n,i-1}(x))/{{\delta}}\quad\mbox{ for }i\in \Lambda.
\end{eqnarray}
Then \eqref{e:BVP_1} and \eqref{e:BVP_i} become
\begin{equation}\label{e:delta2_ui_is_a_uni_sol}
z_{n,i}-\Delta_\Omega u_{n,i}+|w_{n,i}|^{p-1}\cdot w_{n,i}=f_{n,i}\qquad\mbox{ for }i\in\Lambda.
\end{equation}
Let $D_T=[0,T]\times\Omega$, $D_{T,i}:=[t_{i-1}, t_i]\times\Omega$ and $\widetilde{D}_{T,i}:=(t_{i-1}, t_i]\times\Omega$ for $i\in\Lambda$. We construct Rothe's sequence $\{u^{(n)}(t,x)\}$ as below
\begin{equation}\label{e:defi_u(n)}
u^{(n)}(t,x)=u_{n,i-1}(x)+(t-t_i)\cdot w_{n,i}(x)\qquad \mbox{for }(t,x)\in D_{T,i}.
\end{equation}
Also, we define the auxiliary functions
\begin{equation}\label{e:defi_delta_u(n)}
w^{(n)}(t,x)=w_{n,i-1}(x)+(t-t_i)\cdot z_{n,i}(x)\qquad \mbox{for }(t,x)\in D_{T,i},
\end{equation}
and some step functions
\begin{equation}\label{e:defi_delta_overline_u(n)}
\overline{u}^{(n)}(t,x)=
\left\{
\begin{aligned}
&u_{n,i}(x),\qquad&&(t,x)\in\widetilde{D}_{T,i},\\
&g(x),\qquad&&(t,x)\in[-{\delta},0]\times\Omega^\circ,\\
&0,\qquad&&(t,x)\in[-{\delta},0]\times \partial\Omega,
\end{aligned}
\right.
\end{equation}

\begin{equation}\label{e:defi_delta_overline_delta_u(n)}
\overline{w}^{(n)}(t,x)=
\left\{
\begin{aligned}
&w_{n,i}(x),\qquad&&(t,x)\in\widetilde{D}_{T,i},\\
&h(x),\qquad&&(t,x)\in[-{\delta},0]\times\Omega^\circ,\\
&0,\qquad&&(t,x)\in[-{\delta},0]\times \partial\Omega,
\end{aligned}
\right.
\end{equation}
\begin{equation}\label{e:defi_f(n)}
f^{(n)}(t,x)=
\left\{
\begin{aligned}
&f(t_i,x),\qquad\qquad&&(t,x)\in\widetilde{D}_{T,i},\\
&f(0,x),\qquad\qquad&& x\in\Omega^\circ,\\
&0,\qquad\qquad&& t=0, x\in\partial\Omega.
\end{aligned}
\right.
\end{equation}

In order to show that Rothe's sequence $\{u^{(n)}(t,x)\}$ is convergent, more precisely, the sequence converges to $u(t,x)$, a solution of \eqref{e:wave_eq_on_graph}, we give some priori estimates in the following lemma. From now on, we assume that \eqref{eq:condition_on_f} holds.

\begin{lem}\label{L:priori_estimates_eq}
There exist an integer $N_0>0$ and positive constants $C_\Omega$ and $C_{\Omega, p}$ such that for any $n\ge N_0$ and any $i\in \Lambda$,
\begin{eqnarray}\label{e:bounded_of_delta ui_and_ui_and_delta2_ui}
\begin{aligned}
&\|w_{n,i}\|^2_{L^2(\Omega)}+\|\nabla u_{n,i}\|^2_{L^2(\Omega)}+\|u_{n,i}\|^2_{L^2(\Omega)}
+\|w_{n,i}\|^{2}_{L^{2p}(\Omega)}\le C_\Omega,\qquad\|z_{n,i}\|^2_{L^2(\Omega)}\le C_{\Omega,p}.
\end{aligned}
\end{eqnarray}
\end{lem}

\begin{proof}
In view of assumption \eqref{eq:condition_on_f}, we get
\begin{equation*}
\|f(t,\cdot)\|^2_{L^2(\Omega^\circ)}\le C_{\Omega^\circ}T^{2\gamma}+c'\quad\mbox{ for any }t\in[0,T],
\end{equation*}
where $c':=\|f(0,\cdot)\|^2_{L^2(\Omega^\circ)}$.
From \eqref{e:delta2_ui_is_a_uni_sol}, we get for any $i\in \Lambda$ and any $v\in W_0^{1,2}(\Omega)$,
\begin{equation*}
\int_{\Omega^\circ}(z_{n,i}-\Delta_\Omega u_{n,i}+|w_{n,i}|^{p-1}\cdot w_{n,i}-f_{n,i})\cdot v\,d\mu=0.
\end{equation*}
Substituting $v=w_{n,i}$ into the above equation, Lemma \ref{L:Green_formula} implies that
\begin{eqnarray*}
(1-{\delta})\big(\|\nabla u_{n,i}\|^2_{L^2(\Omega)}+\|w_{n,i}\|^2_{L^2(\Omega^\circ)}\big)
\le\|\nabla u_{n,i-1}\|^2_{L^2(\Omega)}+\|w_{n,i-1}\|^2_{L^2(\Omega^\circ)}+{\delta}\|f_{n,i}\|^2_{L^2(\Omega^\circ)}
\end{eqnarray*}
Choosing an integer $N_0>0$ such that ${\delta}<1$ for any $n\ge N_0$, we get
\begin{eqnarray*}
\begin{aligned}
&\|\nabla u_{n,i}\|^2_{L^2(\Omega)}+\|w_{n,i}\|^2_{L^2(\Omega^\circ)}\\
\le&(1-{\delta})^{-i}\Big(\|\nabla u_{n,0}\|^2_{L^2(\Omega)}+\|w_{n,0}\|^2_{L^2(\Omega^\circ)}
+{\delta}\sum_{k=1}^i(1-{\delta})^{k-1}\|f_{n,k}\|^2_{L^2(\Omega^\circ)}\Big)\\
\le&(1-{\delta})^{-n}\Big(\|\nabla u_{n,0}\|^2_{L^2(\Omega)}+\|w_{n,0}\|^2_{L^2(\Omega^\circ)}
+{\delta}\sum_{k=1}^i\|f_{n,k}\|^2_{L^2(\Omega^\circ)}\Big)\\
\le&e^T\Big(\|\nabla u_{n,0}\|^2_{L^2(\Omega)}+\|w_{n,0}\|^2_{L^2(\Omega^\circ)}
+T(C_{\Omega^\circ}T^{2\gamma}+c')\Big)
\le C_\Omega.
\end{aligned}
\end{eqnarray*}

Theorem \ref{T:Sobolev_embedding_thm} implies that $\|u_{n,i}\|^2_{L^2(\Omega^\circ)}\le C_\Omega\|\nabla u_{n,i}\|^2_{L^2(\Omega)}\le C_\Omega^2$. Also,
$$\big(\int_{\Omega}|w_{n,i}|^{2p}\,d\mu\big)^{1/p}\le C_\Omega^{2}\int_{\Omega}|\nabla w_{n,i}|^2\,d\mu\quad\mbox{ for }p>1.$$
Since $\|w_{n,i}\|^2_{L^2(\Omega)}\le C_\Omega$, we have $|w_{n,i}(x)|\le \sqrt{C_\Omega/{\mu_0}}$, and so
\begin{equation*}
\int_{\Omega}|\nabla w_{n,i}|^2\,d\mu\le 4 D_\mu C_\Omega \mu(\Omega)/{\mu_0},
\end{equation*}
where $\mu_0=\min_{x,y\in \Omega}\omega_{xy}$.
This leads to
\begin{equation*}
\|w_{n,i}\|^{2}_{L^{2p}(\Omega)}\le 4 D_\mu C_\Omega^3 \mu(\Omega)/{\mu_0}.
\end{equation*}

The fact $|\Delta_\Omega u_{n,i}(x)|^2\le D_\mu|\nabla u_{n,i}(x)|^2$ implies that
\begin{eqnarray*}
\int_{\Omega^\circ}|\Delta_\Omega u_{n,i}(x)|^2\,d\mu\le C_\Omega D_\mu.
\end{eqnarray*}
It follows from \eqref{e:delta2_ui_is_a_uni_sol} that
\begin{eqnarray*}\label{e:bounded_of_delta2 ui}
\|z_{n,i}\|^2_{L^2(\Omega)}
\le 2\Big(\int_{\Omega^\circ}|\Delta_\Omega u_{n,i}|^2\,d\mu+\int_\Omega|w_{n,i}|^{2p}\,d\mu\Big)
\le C_{\Omega, p}.
\end{eqnarray*}
The proof of Lemma \ref{L:priori_estimates_eq} is completed.
\end{proof}

According to Lemma \ref{L:priori_estimates_eq}, we get the following result.
\begin{lem}
For any $t\in[0,T]$, any $n\ge N_0$ and constants $C_\Omega$, $C_{\Omega, p}$,
\begin{eqnarray}\label{e:bounded_1}
\begin{aligned}
&\|u^{(n)}(t,\cdot)\|_{L^{2}(\Omega)}+\|\overline{u}^{(n)}(t,\cdot)\|_{L^{2}(\Omega)}
+\|w^{(n)}(t,\cdot)\|_{L^2(\Omega)}\\
+&\|\overline{w}^{(n)}(t,\cdot)\|_{L^2(\Omega)}+\|\overline{w}^{(n)}(t,\cdot)\|_{L^{2p}(\Omega)}\le C_\Omega,
\end{aligned}
\end{eqnarray}
\begin{eqnarray}\label{e:bounded_2}
\|w_t^{(n)}(t,\cdot)\|_{L^2(\Omega)}\le C_{\Omega,p}.
\end{eqnarray}
\begin{eqnarray}\label{e:bounded_3}
\|u^{(n)}(t,\cdot)-\overline{u}^{(n)}(t,\cdot)\|_{L^2(\Omega)}\le C_{\Omega}/n
\end{eqnarray}
\begin{eqnarray}\label{e:bounded_4}
\|w^{(n)}(t,\cdot)-\overline{w}^{(n)}(t,\cdot)\|_{L^2(\Omega)}\le C_{\Omega,p}/n.
\end{eqnarray}
\end{lem}

\begin{lem}\label{L:limits_eq}
There exist a function $u\in L^2(\Omega)$ satisfying $u_t, u_{tt}\in L^2(\Omega)$, and two subsequences $\{u^{(n_k)}\}$, $\{\overline{u}^{(n_k)}\}$ such that for any $(t,x)\in D_T$,
\begin{enumerate}
\item[(a)] $u^{(n_k)}\to u$ and $\overline{u}^{(n_k)}\to u$;
\item[(b)] $w^{(n_k)}\to u_t$ and $\overline{w}^{(n_k)}\to u_t$;
\item[(c)] $w_t^{(n_k)}\to u_{tt}$.
\end{enumerate}
\end{lem}
\begin{proof}
(a) Since $\|u^{(n)}\|_{L^2(\Omega)}$ and $\|\overline{u}^{(n)}\|_{L^2(\Omega)}$ are bounded, we have
\begin{eqnarray*}
u^{(n_k)}(t,\cdot)\to u(t,\cdot),\quad\overline{u}^{(n_k)}(t,\cdot)\to \overline{u}(t,\cdot)\qquad\mbox{in }L^2(\Omega)
\end{eqnarray*}
for two subsequences $\{u^{(n_k)}\}, \{\overline{u}^{(n_k)}\}$ and two functions $u, \overline{u}$.
This leads to
\begin{eqnarray}\label{e:unk_point_conv_to_u}
u^{(n_k)}(t,x)\to u(t,x),\qquad\overline{u}^{(n_k)}(t,x)\to \overline{u}(t,x)\quad\mbox{ on }D_T.
\end{eqnarray}
Since $u^{(n_k)}, \overline{u}^{(n_k)}\in W_0^{1,2}(\Omega)$, using \eqref{e:unk_point_conv_to_u}, we have $u=\overline{u}=0$ on $[0,T]\times\partial\Omega$.
It follows from \eqref{e:bounded_3} and \eqref{e:unk_point_conv_to_u} that
\begin{eqnarray*}
 \|u(t,\cdot)-\overline{u}(t,\cdot)\|^2_{L^2(\Omega)}=\lim_{k\to\infty}\|u^{(n_k)}(t,\cdot)-\overline{u}^{(n_k)}(t,\cdot)\|^2_{L^2(\Omega)}=0\qquad
 \mbox{ on }[0,T].
\end{eqnarray*}
Hence $u=\overline{u}$ on $D_T$. This proves (a).

(b) Similar to (a), there exist two subsequences $\{w^{(n_k)}\}$, $\{\overline{w}^{(n_k)}\}$ and a function $w\in L^2(\Omega)$ such that
\begin{eqnarray}\label{e:delta_over_unk_conv}
w^{(n_k)}(t,x)\to w(t,x)\quad\mbox{ and }\quad\overline{w}^{(n_k)}(t,x)\to w(t,x)\quad\mbox{on }D_T.
\end{eqnarray}
Also, $w=0$ on $[0,T]\times\partial\Omega$. Note that for any $t\in[t_{i-1}, t_i]\subseteq[0,T]$ and any $x\in\Omega^\circ$,
\begin{eqnarray*}
\begin{aligned}
u^{(n_k)}(t,x)-g(x)
=&\int_0^{t_1}u_s^{(n_k)}(s,\cdot)\,ds+\cdots+\int_{t_{i-2}}^{t_{i-1}}u_s^{(n_k)}(s,\cdot)\,ds+\int_{t_{i-1}}^tu_s^{(n_k)}(s,\cdot)\,ds\\
=&\int_0^{t_1}w_{n,1}(\cdot)\,ds+\cdots+\int_{t_{i-2}}^{t_{i-1}}w_{n,i-1}(\cdot)\,ds
+\int_{t_{i-1}}^tw_{n,i}(\cdot)\,ds\\
=&\int_0^t\overline{w}^{(n_k)}(s,x)\,ds.
\end{aligned}
\end{eqnarray*}
Letting $k\to\infty$, we get
$$u(t,x)-g(x)=\int_0^t w(s,x)\,ds,$$
where we use
$$\int_0^t\overline{w}^{(n_k)}(s,x)\,ds\to\int_0^t w(s,x)\,ds\quad\mbox{ on }[0,T],$$
which follows from $\overline{w}^{(n_k)}$ is bounded on $D_T$ and Dominated Convergence Theorem.
Hence $w=u_t$, $u(0,x)=g(x)$ for $x\in\Omega^\circ$ and $u_t=0$ on $[0,T]\times\partial\Omega$.

(c) Similar to (a), there exist a subsequence $\{w_t^{(n_k)}\}$ satisfying
\begin{equation*}
w_t^{(n_k)}(t,\cdot)\to u_{tt}\qquad\mbox{ on }D_T.
\end{equation*}
Also, $u_t|_{t=0}=h$ on $\Omega^\circ$. In the proof, we use the fact that
\begin{eqnarray}\label{e:limit_w_t_to_limit_u_tt}
\int_0^t w_s^{(n_k)}(s,x)\,ds\to \int_0^t u_{ss}^{(n_k)}(s,x)\,ds\quad\mbox{on }D_T.
\end{eqnarray}
\end{proof}

\begin{lem}\label{L:int_limits_eq}
The following results hold:
\begin{enumerate}
\item[(a)] $\int_0^T \Delta_\Omega\overline{u}^{(n_k)}(t,x)\,dt\to\int_0^T\Delta_\Omega u(t,x)\,dt$ on $\Omega^\circ$;
\item[(b)] $\int_0^T |\overline{w}^{(n_k)}(t,x)|^{p-1}\cdot \overline{w}^{(n_k)}(t,x)\,dt\to\int_0^T|u_t(t,x)|^{p-1}\cdot u_t(t,x)\,dt$ on $\Omega$;
\item[(c)] $\int_0^T f^{(n_k)}(t,x)\,dt\to\int_0^T f(t,x)\,dt$ on $\Omega^\circ$.
\end{enumerate}
\end{lem}
\begin{proof}
(a) It follows from \eqref{e:unk_point_conv_to_u} that $\Delta_\Omega \overline{u}^{(n_k)}(t,x)\to \Delta_\Omega u(t,x)$ on $[0,T]\times\Omega^\circ$. In view of \eqref{e:bounded_1}, we get $\Delta_\Omega \overline{u}^{(n_k)}$ is bounded on $[0,T]\times\Omega^\circ$. Dominated Convergence Theorem implies that (a) holds.

(b,c) The proofs are the same as that of (a).
\end{proof}

\subsection{Proof of Theorem~\ref{T:sol_of_wave_equ}}

Using notation and results in subsection \ref{SS:some_priori_est}, we prove our main theorem.

\begin{proof}[Proof of Theorem~\ref{T:sol_of_wave_equ}]
\noindent{\bf Existence}

In view of \eqref{e:delta2_ui_is_a_uni_sol}, we get for $p>1$,
\begin{equation*}
\int_0^T(z_{n,i}-\Delta_\Omega u_{n,i}+|w_{n,i}|^{p-1}\cdot w_{n,i}-f_{n,i})\,dt=0\qquad\mbox{on }\Omega^\circ.
\end{equation*}
Combining this with \eqref{e:defi_u(n)}--\eqref{e:defi_delta_overline_delta_u(n)}, we obtain
\begin{equation*}
\int_0^T\big(w_s^{(n)}(t,x)-\Delta_\Omega\overline{u}^{(n)}(t,x)
+|\overline{w}^{(n)}(t,x)|^{p-1}\cdot\overline{w}^{(n)}(t,x)-f^{(n)}(t,x)\big)\,dt=0\qquad\mbox{on }\Omega^\circ.
\end{equation*}
Let $u$ be the limit function in Lemma \ref{L:limits_eq}.
Letting $n=n_k$ and taking the limits as $k\to\infty$ in the above equation, Lemma \ref{L:int_limits_eq} and \eqref{e:limit_w_t_to_limit_u_tt} imply that
\begin{equation*}
\int_0^T\big(u_{tt}(t,x)-\Delta_\Omega u(t,x)+|u_t(t,x)|^{p-1}\cdot u_t(t,x)-f(t,x)\big)\,dt=0.
\end{equation*}
From Lemma \ref{L:limits_eq}, we get the initial and boundary conditions of \eqref{e:wave_eq_on_graph} hold. $u$ is a solution of \eqref{e:wave_eq_on_graph} follows from the arbitrary of $T$.

\noindent{\bf Uniqueness}

Let $u$ and $\check{u}$ be two solution of \eqref{e:wave_eq_on_graph}. Let $\varphi:=u-\check{u}$. Then for $p>1$,
\begin{equation*}
\left\{
\begin{aligned}
&\varphi_{tt}-\Delta_\Omega \varphi+|u_t|^{p-1}\cdot u_t-|\check{u}_t|^{p-1}\cdot \check{u}_t=0,\quad&&t\ge0, x\in\Omega^\circ,\\
&\varphi|_{t=0}=0,\quad&& \Omega^\circ,\\
&\varphi_t|_{t=0}=0,\quad&&\Omega^\circ,\\
&\varphi=0,\quad&&t\ge 0, x\in\partial\Omega^\circ.
\end{aligned}
\right.
\end{equation*}
For $t\in[0,\infty)$, let $$G(t):=\int_{\Omega}|\nabla \varphi(t,x)|^2\,d\mu+\int_{\Omega^\circ}|\varphi_t(t,x)|^2\,d\mu.$$
Then $G(0)=0$. Moreover,
\begin{eqnarray*}
\begin{aligned}
G'(t)
=&2\int_{\Omega} \Gamma(\varphi,\varphi_t)\,d\mu
+2\int_{\Omega^\circ}\varphi_t\cdot
\big[\Delta_\Omega \varphi-\big(|u_t|^{p-1}\cdot u_t-|\check{u}_t|^{p-1}\cdot\check{u}_t\big)\big]\,d\mu\\
=&-2\int_{\Omega^\circ} (u_t-\check{u}_t)\cdot\big(|u_t|^{p-1}\cdot u_t-|\check{u}_t|^{p-1}\cdot\check{u}_t\big)\,d\mu\\
\le&0,
\end{aligned}
\end{eqnarray*}
where we use the fact that for $p>1$, $(u_t-\check{u}_t)\cdot\big(|u_t|^{p-1}\cdot u_t-|\check{u}_t|^{p-1}\cdot\check{u}_t\big)\ge 0$.

For any $t\ge 0$, $G'(t)\le 0$ and $G(0)=0$ imply that $G(t)\equiv0$, and hence
$$\nabla \varphi\equiv 0\quad\mbox{on }[0,\infty)\times\Omega\qquad\mbox{and}
\qquad \varphi_t\equiv 0\quad\mbox{on }[0,\infty)\times\Omega^\circ,$$
which together with $\varphi(t,x)=0$ for $t\ge0$ and $x\in\partial\Omega$ and $\varphi(0,x)=0$ for $x\in\Omega^\circ$, we have
$\varphi\equiv 0$. Then $u\equiv\check{u}$ follows.
\end{proof}

\textbf{Acknowledgement}
This work is supported by the National Science Foundation of China [12071245].

Yong Lin,\\
Yau Mathematical Sciences Center, Tsinghua University\\ Beijing, 100084, China.\\
\textsf{E-mail: yonglin@tsinghua.edu.cn}\\
\\
Yuanyuan Xie,\\
School of Mathematics, Renmin University of China\\ Beijing, 100872, China.\\
\textsf{E-mail: yyxiemath@163.com}\\

\end{document}